\documentclass[12pt]{amsart}

\title[Root Systems and Quantum cohomology]{Root Systems and the
Quantum Cohomology of ADE resolutions} 

\author{Jim Bryan and Amin Gholampour}
\address{
Department of Mathematics\\
University of British Columbia \\
Room 121, 1984 Mathematics Road  \\
Vancouver, B.C., Canada V6T 1Z2  
}

\usepackage{diagrams}
\usepackage{verbatim}
\usepackage{amsmath}
\usepackage{amsmath,amsthm,amsfonts}

\newcommand{\cnums} {{\mathbb C}}          
\newcommand{\znums} {{\mathbb Z}}		
\newcommand{\qnums} {{\mathbb Q}}		

\newcommand{\tr}{\operatorname{tr}}

\newcommand{\Def}{\operatorname{Def}}

\renewcommand{\P}{\mathbb{P}}
\newcommand{\M}{\overline{{M}}}

\newcommand{\LangRang}[1]{\left\langle #1 \right\rangle}
\newcommand{\abar}{\overline{\alpha }}
\newcommand{\bbar}{\overline{\beta }}
\newcommand{\gbar}{\overline{\gamma }}
\newcommand{\qbar}{\overline{q}}
\renewcommand{\tilde}{\widetilde}
\renewcommand{\O}{\mathcal{O}}
\newcommand{\X}{\mathcal{X}}
\newcommand{\orb}{\mathrm{orb}}
\newcommand{\stab}{\operatorname{stab}}

\newtheorem{thm}{Theorem}
\newtheorem{theorem}[thm]{Theorem}
\newtheorem{cor}[thm]{Corollary}

\newtheorem{lemma}[thm]{Lemma}

\newtheorem{proposition}[thm]{Proposition}

\newtheorem{defn}[thm]{Definition}

\newtheorem{conjecture}[thm]{Conjecture}
\newtheorem{conjecture*}{Conjecture}

\begin{document}

\begin{abstract}
We compute the $\cnums^*$-equivariant quantum cohomology ring of $Y$,
the minimal resolution of the DuVal singularity $\cnums^2/G$ where $G$
is a finite subgroup of $SU(2)$. The quantum product is
expressed in terms of an ADE root system canonically associated to
$G$. We generalize the resulting Frobenius manifold to non-simply
laced root systems to obtain an $n$ parameter family of algebra
structures on the affine root lattice of any root system. Using the
Crepant Resolution Conjecture, we obtain a prediction for the orbifold
Gromov-Witten potential of $[\cnums ^{2}/G]$.
\end{abstract}
\maketitle 



\section{Introduction}\label{sec: intro}

\subsection{Overview.}
Let $G$ be a finite
subgroup of $SU(2)$, and let
\[
Y\to \cnums ^{2}/G
\]
be the minimal resolution of the corresponding DuVal singularity.  The
classical McKay correspondence describes the geometry of $Y$ in terms
of the representation theory of $G$
\cite{Gonzalez-Sprinberg-Verdier,McKay,Reid-asterisque}.

The geometry of $Y$ gives rise to a Dynkin diagram of ADE type. The
nodes of the diagram correspond to the irreducible components of the
exceptional divisor of $Y $. Two nodes have a connecting edge if
and only if the corresponding curves intersect.

Associated to every Dynkin diagram of ADE type is a simply laced root
system. In this paper, we describe the $\cnums^*$-equivariant quantum
cohomology of $Y$ in terms of the associated root system. This
provides a quantum version of the classical McKay correspondence.

\subsection{Results.}

The set $\{E_1,\ldots,E_n\}$ of irreducible components of the
exceptional divisor of $Y$ forms a basis of $H_{2} (Y,\znums )$. The
intersection matrix $E_{i}\cdot E_{j}$ defines a perfect pairing on
$H_{2} (Y,\znums )$. Let $R$ be the simply laced root system
associated to the Dynkin diagram of $Y$. We can identify $H_{2}
(Y,\znums )$ with the root lattice of $R$ in a way so that
$E_{1},\dotsc ,E_{n}$ correspond to simple roots $\alpha
_{1},\dotsc ,\alpha _{n}$ and the intersection matrix is minus the
Cartan matrix
\[
E_{i}\cdot E_{j} = -\LangRang{\alpha _{i},\alpha _{j}}.
\]

Using the above pairing, we identify $H^{2} (Y,\znums )$ with $H_{2}
(Y,\znums )$ (and hence with the root lattice). Since the scalar
action of $\cnums ^{*}$ on $\cnums ^{2}$ commutes with the action of $G$,
$\cnums ^{*}$ acts on $\cnums ^{2}/G$ and this action lifts to an action on
$Y$. The cycles $E_{1},\dotsc ,E_{n}$ are $\cnums ^{*}$
invariant, and so the classes $\alpha _{1},\dotsc ,\alpha _{n}$ have natural
lifts to equivariant (co)homology.  Additively, the equivariant
quantum cohomology ring is thus a free module generated by the classes
$\{1,\alpha _{1},\dotsc ,\alpha _{n} \}$. The ground ring is
$\znums[t][[q_1,\ldots,q_n]]$ where $t$ is the equivariant parameter
and $q_{1},\dotsc ,q_{n}$ are the quantum parameters associated to the
curves $E_{1},\dotsc ,E_{n}$. So additively we have
\[
QH_{\cnums^*}^*(Y)\cong
H^*(Y,\znums)\otimes
\znums[t][[q_1,\ldots,q_n]].
\]
We extend the pairing $\LangRang{\,,\,}$ to a $\qnums
[t,t^{-1}][[q_{1},\dotsc ,a_{n}]]$ valued pairing on $QH_{\cnums ^{*}}^{*} (Y)$ by
making 1 orthogonal to $\alpha _{i}$ and setting
\[
\LangRang{1,1} = \frac{-1}{t^{2}|G|}.
\]
The product structure of $QH^{*}_{\cnums ^{*}} (Y)$ is determined by
our main theorem:

\bigskip

\begin{theorem}\label{thm: product in QH}
Let $v,w\in H^2(Y,\znums)$ which we identify with the root lattice of
$R$ as above. Then the quantum product of $v$ and $w$ is given by the
formula:
\[
v\star w = -t^{2}|G|\LangRang{v,w} +\sum_{\beta\in R^+} \LangRang{v ,\beta}\LangRang{w,\beta }
t\frac{1+q^\beta}{1-q^\beta}\beta
\]
where the sum is over the positive roots of $R$ and for $\beta =\sum
_{i=1}^{n}b_{i}\alpha _{i}$, $q^{\beta }$ is defined by
\[
q^\beta=\prod _{i=1}^{n}q_{i}^{b_{i}}
\]
The quantum product satisfies the Frobenius condition
\[
\LangRang{v\star w,u} =\LangRang{v,w\star u}
\]
making $QH_{\cnums ^{*}}^{*} (Y)$ a Frobenius algebra over $\qnums
[t,t^{-1}][[q_{1},\dotsc ,q_{n}]]$.
\end{theorem}

\bigskip

Note that by a standard fact in root theory \cite[VI.1.1 Proposition~3
and V.6.2 Corollary to Theorem~1]{Bourbaki}, the formula in
Theorem~\ref{thm: product in QH} can alternatively be written as
\[
v\star w = \sum_{\beta\in R^+} \LangRang{v ,\beta}\LangRang{w,\beta }
\left(-t^2\frac{|G|}{h}+t\frac{1+q^\beta}{1-q^\beta}\beta\right)
\]
where $h=\frac{|R|}{n}$ is the Coxeter number of $R$.

We remark that we can regard $H^0(Y)\oplus H^2(Y)$ as the root lattice
for the affine root system and consequently, we can regard
$QH_{\cnums^*}^*(Y)$ as defining a family of algebra structures on the
affine root lattice depending on variables $t,q_1,\ldots,q_n$. We also
remark that even though the product in Theorem~\ref{thm: product in
QH} is expressed purely in terms of the root system, we know of no
root theoretic proof of associativity, even in the ``classical'' limit
$q_{i}\to 0$. 

In section~\ref{sec: non-simply laced case}, which can be read
independently from the rest of this paper, we will generalize our
family of algebras to root systems which are not simply-laced
(Theorem~\ref{thm: QH for any R}). We will prove associativity of the
product in the non-simply laced case by reducing it to the simply
laced case. Our formula also allows us to prove that the action of the
Weyl group induces automorphisms of the Frobenius algebra
(Corollary~\ref{cor: W acts on QH by automorphism}).

Our theorem is formulated as computing \emph{small} quantum
cohomology, but since the cohomology of $Y$ is concentrated in degree
0 and degree 2, the large and small quantum cohomology rings contain
equivalent information. The proof of Theorem~\ref{thm: product in QH}
requires the computations of genus 0 equivariant Gromov-Witten
invariants of $Y$. This is done in section~\ref{sec: GW theory of Y}.

In section~\ref{sec: CRC}, we use the Crepant Resolution Conjecture
\cite{Bryan-Graber} and our computation of the Gromov-Witten
invariants of $Y$, to obtain a prediction for the orbifold
Gromov-Witten potential of $[\cnums ^{2}/G]$ (Conjecture~\ref{conj:
prediction for F_X}).

\subsection{Relationship to other work} A certain specialization of
the Frobenius algebra $QH_{\cnums ^{*}}^{*} (Y)$ appears as the
quantum cohomology of the $G$-Hilbert scheme resolution of $\cnums
^{3}/G$ for $G\subset SO (3)$, (see \cite{Bryan-Gholampour3}). The
equivariant Gromov-Witten theory of $Y$ in higher genus has been
determined by recent work of Maulik \cite{Maulik-An}.

\section{Gromov-Witten theory of $Y$}\label{sec: GW theory of Y}

In this section we compute the equivariant genus zero Gromov-Witten
invariants of $Y$. The invariants of non-zero degree are computed by
relating them to the invariants of a certain threefold $W$ constructed
as the total space of a family of deformations of $Y$. The invariants
of $W$ are computed by the method of Bryan, Katz, and Leung
\cite{BKL}. The degree zero invariants are computed by localization.

\subsection{Invariants of non-zero degree.}

$\Def(Y)$, the versal space of $\cnums^*$-equivariant deformations of
$Y$ is naturally identified with the complexified root space of the
root system $R$ \cite{Ka-Mo}. A generic deformation of $Y$ is an
affine variety and consequently has no compact curves. The hyperplane
$D_\beta\subset \Def(Y)$ perpendicular to a positive root 
\[
\beta = \sum_{i=1}^n b_i\alpha_i
\]
parameterizes those deformations of $Y$ for which the curve
\[
b_1E_1+\cdots + b_nE_n
\]
also deforms. Moreover, for a generic point $t\in D_{\beta }$, the
corresponding curve is a smooth $\P ^{1}$ which generates the Picard
group of the corresponding surface (\cite[Prop.~2.2]{BKL} and
\cite[Thm.~1]{Ka-Mo}).

Let 
\[
\imath:\cnums \to \Def(Y)
\]
be a generic linear subspace. We obtain a threefold $W$ by pulling
back the universal family over $\Def(Y)$ by $\imath$. The embedding
$\imath $ can be made $\cnums^*$-equivariant by defining the action on
$\cnums$ to have weight $2$. This follows from \cite[Thm.~1]{Ka-Mo}
after noting that the $\cnums ^{*} $ action in \cite{Ka-Mo} is the
square of the action induced by the action on $\cnums ^{2}/G$. Clearly
$Y\subset W$ and the degree of the normal bundle is
\[
c_1(N_{Y/W})=2t
\]
(recall that $t$ is the equivariant parameter).

The threefold $W$ is Calabi-Yau and its Gromov-Witten invariants are
well defined in the non-equivariant limit. This assertion follows from
the fact that the moduli space of stable maps to $W$ is compact. This
in turn follows from the fact that $W$ admits a birational map $W\to
W_{\mathrm{aff}}$ contracting $E_{1}\cup \dotsb \cup E_{n}$ such that
$W_{\mathrm{aff}}$ is an affine variety (see
\cite{BKL,Ka-Mo}). Consequently, all non-constant stable maps to $W$
must have image contained in the exceptional set of $W\to
W_{\mathrm{aff}}$ and thus, in particular, all non-constant stable
maps to $W$ have their image contained in $Y$.

There is a standard technique in Gromov-Witten theory for comparing
the virtual class for stable maps to a submanifold to the virtual
class for the stable maps to the ambient manifold when all the maps
have image contained in the submanifold \cite{Be-Fa}.  This allows us
to compare the Gromov-Witten invariants of $W$ and $Y$.

For any non-zero class 
\[
A\in H_{2} (Y)\subset H_{2} (W)
\]
let
\[
\LangRang{\,\,}_{A}^{Y},\quad \LangRang{\,\,}_{A}^{W}
\]
denote the genus zero, degree $A$, zero insertion Gromov-Witten
invariant of $Y$ and $W$ respectively. We have
\[
\LangRang{\,\,}_{A}^{W}= \int _{[\M _{0,0}
(Y,A)]^{vir}}e (-R^{\bullet }\pi _{*}f^{*}N_{Y/W})
\]
where $\M _{0,0} (Y,A)$ is the moduli space of stable maps, $\pi :C\to
\M _{0,0} (Y,A)$ is the universal curve, $f:C\to Y$ is the universal
map, and $e$ is the equivariant Euler class.

Since the line bundle $N_{Y/W} $ is trivial up to the $\cnums ^{*}$
action, and $\pi $ is a family of genus zero curves, we get
\[
R^{\bullet }\pi _{*}f^{*}N_{Y/W} = R^{0}\pi _{*}f^{*}N_{Y/W}= 
\O\otimes \cnums _{2t}
\]
where $\cnums _{2t}$ is the $\cnums ^{*}$ representation of weight 2
so that we have
\[
c_{1} (\O\otimes \cnums _{2t}) = 2t.
\]
Consequently, we have
\[
e (-R^{\bullet }\pi _{*}f^{*}N_{Y/W}) = \frac{1}{2t}
\]
and so
\begin{align*}
\LangRang{\,\,}_{A}^{W}&=\int _{[\M _{0,0} (Y,A)]^{vir}}\frac{1}{2t}\\
&=\frac{1}{2t}\LangRang{\,\,}_{A}^{Y}.
\end{align*}

To compute $\LangRang{\,\,}_{A}^{W}$, we use the deformation
invariance of Gromov-Witten invariants. Although $W$ is non-compact,
the moduli space of stable maps is compact, and the deformation of $W$
is done so that the stable map moduli spaces are compact throughout
the deformation. The technique is identical to the deformation
argument used in \cite{BKL} where it is presented in greater detail.

We deform $W$ to a threefold $W'$ as follows. Let
\[
\imath ' :\cnums \to \Def (Y)
\]
be a generic affine linear embedding and let $W'$ be the pullback by
$\imath '$ of the universal family over $\Def (Y)$. The threefold $W'$
is a deformation of $W$ since $\imath '$ is a deformation of $\imath
$.

\begin{lemma}
The compact curves of $W'$ consist of isolated $\P ^{1}$s, each having
normal bundle $\O (-1)\oplus \O (-1)$, one in each homology class
$\beta \in H_{2} (W')\cong H_{2} (Y)$ corresponding to a positive root.
\end{lemma}
\textsc{Proof:} The map $\imath '$ intersects each hyperplane
$D_{\beta }$ transversely in a single generic point $t$. The surface
$S_{t}$ over the point $t$ contains a single curve $C_{t}\cong \P
^{1}$ of normal bundle $N_{C_{t}/S_{t}}\cong \O (-2)$
and this curve is in the class $\beta $. There is a short exact
sequence
\[
0\to N_{C_{t}/S_{t}}\to  N_{C_{t}/W'}\to \O \to 0
\]
and since $\imath '$ intersects $D_{\beta }$ transversely, $C_{\beta
}$ does not have any deformations (even infinitesimally) inside
$W'$. Consequently, we must have $N_{C_{\beta }/W'}\cong \O (-1)\oplus
\O (-1)$.\qed 

Since all the curves in $W'$ are isolated $(-1,-1)$ curves, we can
compute the Gromov-Witten invariants of $W'$ using the
Aspinwall-Morrison multiple cover formula. Combined with the
deformation invariance of Gromov-Witten invariants, we obtain:

\begin{lemma}\label{lem: invariants of W}
For $A\neq 0$ we have
\[
\LangRang{\,\,}_{A}^{Y} = 2t\LangRang{\,\,}_{A}^{W}=2t\LangRang{\,\,}_{A}^{W'}=
\begin{cases}
2t\frac{1}{d^{3}}&\text{ if $A=d\beta $ where $\beta$ is a positive root}\\
0&\text{otherwise.}
\end{cases}
\]
\end{lemma}

Since all the cohomology of $Y$ is in $H^{0} (Y)$ and $H^{2} (Y)$, the
$n$-point Gromov-Witten invariants of non-zero degree are determined
from the 0-point invariants by the divisor and the fundamental class
axioms.

\subsection{Degree 0 invariants.}
The only non-trivial degree zero invariants have 3 insertions and are
determined by classical integrals on $Y$. They are given in the
following lemma.

\begin{lemma}\label{lem: classical invs}
Let 1 be the generator of $H^{0}_{\cnums ^{*}} (Y)$ and let $\{\alpha
_{1},\dotsc ,\alpha _{n} \}$ be the basis for $H^{2}_{\cnums ^{*}}
(Y)$ which is also identified with the simple roots of $R$ as in
section~\ref{sec: intro}. Then the degree 0, 3-point Gromov-Witten
invariants of $Y$ are given as follows:
\begin{align}
\LangRang{1,1,1}_{0} &=\frac{1}{t^{2}|G|}\label{eqn: <111>},\\
\LangRang{\alpha _{i},1,1}_{0}&=0\label{eqn: <11ai>},\\
\LangRang{\alpha _{i},\alpha _{j},1}_{0}&=-\LangRang{\alpha _{i},\alpha _{j}} \label{eqn: <1aiaj>},\\
\LangRang{\alpha _{i},\alpha _{j},\alpha _{k}}_{0} &=-t\sum _{\beta \in
R^{+}} \LangRang{\alpha _{i},\beta } \LangRang{\alpha _{j},\beta }
\LangRang{\alpha _{k},\beta }.\label{eqn: <aiajak>}
\end{align}
\end{lemma}

\textsc{Proof:} The degree zero, genus zero, 3-point Gromov-Witten
invariants are given by integrals over $Y$:
\[
\LangRang{x, y, z}_{0} =\int _{Y}x\cup y\cup z.
\]
Because $Y$ is non-compact, the integral must be defined\footnote{We
remark that this method of defining the Gromov-Witten invariants of a
non-compact space does not affect the desired properties of quantum
cohomology: the associativity still holds and the Frobenius structure
still exists with the novelty that the pairing takes values in the
ring $\qnums [t,t^{-1}]$. See \cite[section~1.4]{Bryan-Graber}, for a
discussion.} via $\cnums ^{*}$ localization and takes values in
$\qnums [t,t^{-1}]$, the localized equivariant cohomology ring of a
point:
\begin{diagram}
\int _{Y}\quad : H^{*}_{\cnums ^{*}} (Y) &\rTo &\qnums [t,t^{-1}]\\
\phi  &\mapsto &\int _{F}\frac{\phi |_{F}}{e (N_{F/Y})}.
\end{diagram}
Here $F\subset Y$ is the (compact) fixed point locus of the action of
$\cnums ^{*}$ on $Y$.

By correspondence of residues \cite{Bertram}, integrals over $Y$ can be
computed by first pushing forward to $\cnums ^{2}/G $ followed by (orbifold)
localization on $\cnums ^{2}/G$. Equation~(\ref{eqn: <111>}) follows immediately:
\[
\int _{Y}1 = \int _{\cnums ^{2}/G}1 = \frac{1}{t^{2}|G|}.
\]
The factor $t^{2}$ is the equivariant Euler class of the normal bundle
of $[0/G]\subset [\cnums ^{2}/G]$ and the factor $\frac{1}{|G|}$
accounts for the automorphisms of the point $[0/G]$.

Let $L_{i}\to  Y$ be the $\cnums ^{*}$ equivariant line bundle with 
\[
c_{1} (L_{i})=\alpha _{i}.
\]
Since $\alpha _{i}$ was defined to be dual to $E_{i}$ via the
intersection pairing, we have
\[
\int _{E_{j}}c_{1} (L_{i}) = E_{i}\cdot E_{j} = -\LangRang{\alpha
_{i},\alpha _{j}}.
\]
Computing the left hand side using localization, we see that the
weight of the $\cnums ^{*}$ action on $L_{i}$ at a fixed point $p\in
E_{i}$ must be the same as the weight of the $\cnums ^{*}$ action on
the normal bundle $N_{E_{i}/Y} $ at $p$, and the weight of the action
on $L_{i}$ is 0 over fixed points not on $E_{i}$.

Equation~(\ref{eqn: <1aiaj>}) and Equation~(\ref{eqn: <11ai>}) then
easily follow from localization.

To prove Equation~(\ref{eqn: <aiajak>}), we compute the left hand side
by localization to get
\begin{equation*}
\LangRang{\alpha _{i},\alpha _{j},\alpha _{k}}_{0} = \begin{cases}
0&\text{if }E_{i}\cup E_{j}\cup E_{k} = \emptyset, \\
-8t&\text{if }i=j=k,\\
w_{ijj}&\text{if $i\neq j=k$ and $E_{i}\cup E_{j}\neq \emptyset $}
\end{cases}
\end{equation*}
where 
\[
w_{ijj} =c_{1} (N_{E_{j}/Y}|_{p_{ij}})
\]
is the weight of the $\cnums ^{*}$ action on the normal bundle of
$E_{j}$ at the point $p_{ij}=E_{i}\cup E_{j}$. 

The normal weights $w_{iij}$ satisfy the following three conditions:
\begin{enumerate}
\item Since $K_{Y}$ is the trivial bundle with a $\cnums ^{*}$ action
of weight $2t$, the sum of the normal weights at $p=E_{i}\cap E_{j}$
is $2t$ and so
\[
w_{ijj}+w_{jii}=2t \text{ when $E_{i}\cap E_{j}\neq \emptyset $ and
$i\neq j$.}
\]
\item Since $E_{i}$ is $\cnums ^{*}$ invariant, the sum of the tangent
weights of any two distinct fixed points on $E_{i}$ is zero. Combined
with the above, we see that the sum of the normal weights at any two
distinct fixed points is $4t$ so
\[
w_{ikk} + w_{jkk} = 4t\text{ when $E_{i}\cap E_{k}\neq \emptyset $,
$E_{j}\cap E_{k}\neq \emptyset $, and $i\neq j\neq k$.}
\]
\item Since automorphisms of the Dynkin diagrams induce equivariant
automorphisms of $Y$, the normal weights are invariant under such
automorphisms.
\end{enumerate}

The normal weights are completely determined by the above three
conditions. Indeed, it is clear that once one normal weight is known,
then properties (1) and (2) determine the rest. Moreover, in the case
of Dynkin diagrams of type $D_{n}$ or $E_{n}$, the curve corresponding
to the trivalent vertex of the Dynkin graph must be fixed and so its
tangent weights are zero. In the $A_{n}$ case, condition (3) provides
the needed extra equation.

To summarize the above, the three point degree zero invariants
$\LangRang{\alpha _{i},\alpha _{j},\alpha _{k}}_{0}$ satisfy the
following conditions and are completely determined by them.
\begin{enumerate}
\item [(i)]$\LangRang{\alpha _{i},\alpha _{j},\alpha _{k}}_{0} $ is
symmetric in $\{i,j,k \}$,
\item [(ii)]$\LangRang{\alpha _{i},\alpha _{j},\alpha _{k}}_{0}$ is
invariant under any permutation of indices induced by a Dynkin diagram
automorphism,
\item [(iii)]$\LangRang{\alpha _{i},\alpha _{j},\alpha _{k}}_{0}=0$ if
$\LangRang{\alpha _{j},\alpha _{k}}=0$,
\item [(iv)]$\LangRang{\alpha _{i},\alpha _{j},\alpha _{k}}_{0}=-8t$,
if $i=j=k$,
\item [(v)]$\LangRang{\alpha _{i},\alpha _{i},\alpha
_{j}}_{0}+\LangRang{\alpha _{j},\alpha _{j},\alpha _{i}}_{0}=2t$ if
$\LangRang{\alpha _{i},\alpha _{j}}=-1$,
\item [(vi)]$\LangRang{\alpha _{i},\alpha _{k},\alpha
_{k}}_{0}+\LangRang{\alpha _{j},\alpha _{k},\alpha _{k}}_{0}=4t$ if
$i\neq j$ and $\LangRang{\alpha _{i},\alpha _{k}}=\LangRang{\alpha
_{j},\alpha _{k}}=-1$.
\end{enumerate}

So to finish the proof of Lemma~\ref{lem: classical invs}, it suffices
to show that the right hand side of equation~\eqref{eqn: <aiajak>}
also satisfies all the above properties.  This is precisely the
content of Proposition~\ref{prop: properties of Gijk}, a root
theoretic result which we prove in section~\ref{sec: non-simply laced
case}. \qed

\section{Proof of the main theorem}\label{sec: pf of main thm}

Having computed all the Gromov-Witten invariants of $Y$, we can
proceed to compute the quantum product and prove our main theorem. 

The quantum product $\star $ is defined in terms of the genus 0, 3-point
invariants of $Y$ by the formula:
\[
-\LangRang{x\star y,z} = \sum _{A \in H_{2} (Y,\znums )} \LangRang{x,y,z}_{A }q^{A }
\]
where the strange looking minus sign is due to the fact that the
pairing $\LangRang{\,,\,}$, which coincides with the Cartan pairing on
the roots, is the negative of the cohomological pairing.

To prove our formula for $v\star w$, it suffices to check that the
formula holds after pairing both sides with 1 and with any $u\in H^{2}
(Y)$.

By definition and Lemma~\ref{lem: classical invs} we have
\begin{align*}
-\LangRang{v\star w,1}&=\sum _{A\in H_{2} (Y)} \LangRang{v,w,1}_{A}q^{A}\\
&=\LangRang{v,w,1}_{0}\\
&=-\LangRang{v,w}
\end{align*}
which is in agreement with the right hand side of the formula in
Theorem~\ref{thm: product in QH} when paired with 1 since 1 is
orthogonal to $H^{2} (Y)$ and
\[
\LangRang{1,1} = -\frac{1}{t^{2}|G|}.
\]

For $u\in H^{2} (Y)$ we apply the divisor axiom to get
\begin{align*}
-\LangRang{v\star w,u} &= \sum _{A\in H_{2} (Y)} \LangRang{v,w,u}_{A}q^{A}\\
&=\LangRang{v,w,u}_{0}-\sum _{A\neq
0}\LangRang{v,A}\LangRang{w,A}\LangRang{u,A}\LangRang{\,\, 
}_{A}q^{A}.
\end{align*}
Applying Lemma~\ref{lem: classical invs} and Lemma~\ref{lem:
invariants of W} we get
\begin{align*}
-\LangRang{v\star w,u} &= -t \sum _{\beta \in R^{+}}\LangRang{v,\beta
}\LangRang{w,\beta }\LangRang{u,\beta }-\sum _{\beta \in R^{+}}\sum
_{d=1}^{\infty } \LangRang{v,d\beta } \LangRang{w,d\beta }
\LangRang{u,d\beta }\frac{2t}{d^{3}}q^{d\beta }\\
&= -t\sum _{\beta \in R^{+}}\LangRang{v,\beta } \LangRang{w,\beta }
\LangRang{u,\beta } \left(1+\frac{2q^{\beta }}{1-q^{\beta }} \right)\\
&=-t\sum _{\beta \in R^{+}}\LangRang{v,\beta } \LangRang{w,\beta }
\LangRang{u,\beta }\left(\frac{1+q^{\beta }}{1-q^{\beta }} \right).
\end{align*}
Pairing the right hand side of the formula in Theorem~\ref{thm:
product in QH} with $u$, we find agreement with the above and the
formula for $\star $ is proved.

To prove that the Frobenius condition holds, we only need to observe
that the pairing on $QH_{\cnums ^{*}}^{*} (Y)$ is induced by the three
point invariant with one insertion of 1:
\[
-\LangRang{x,y} = \LangRang{x,y,1}_{0}.
\]
This indeed follows from equations~(\ref{eqn: <111>}), (\ref{eqn:
<11ai>}), and (\ref{eqn: <1aiaj>}).  \qed

\section{The algebra for arbitrary root systems}\label{sec: non-simply
laced case} In this section we construct a Frobenius algebra $QH_{R}$
associated to any irreducible, reduced root system $R$
(Theorem~\ref{thm: QH for any R}). This section can be read
independently from the rest of the paper.

\subsection{Root system notation} In this section we let $R$ be an
irreducible, reduced, rank $n$ root system. That is,
\[
R=\{R,V,\LangRang{\,,\,} \}
\]
consists of a finite subset $R$ of a real inner product space $V$ of
dimension $n$ satisfying
\begin{enumerate}
\item $R$ spans $V$,
\item if $\alpha \in R$ then $k\alpha \in R$ implies $k=\pm 1$,
\item for all $\alpha \in R$, the reflection $s_{\alpha }$ about
$\alpha ^{\perp}$, the hyperplane perpendicular to $\alpha $ leaves
$R$ invariant,
\item for any $\alpha ,\beta \in R$, the number
$\frac{2\LangRang{\alpha ,\beta }}{\LangRang{\alpha ,\alpha }}$ is an
integer, and
\item $V$ is irreducible as a representation of $W$, the Weyl group
(i.e. the group generated by the reflections $s_{\alpha }$, $\alpha \in R$).
\end{enumerate}
We will also assume that the inner product $\LangRang{\,,\,}$ takes
values in $\znums $ on $R$.

Let $\{\alpha _{1},\dotsc ,\alpha _{n} \}$ be a system of simple
roots, namely a subset of $R$ spanning $V$ and such that for every
$\beta =\sum _{i=1}^{n}b_{i}\alpha _{i}$ in $R$ the coefficients
$b_{i}$ are either all non-negative or all non-positive. As is
customary, we define
\[
\alpha ^{\vee }=\frac{2\alpha }{\LangRang{\alpha ,\alpha }}.
\]

We will also require a certain constant $\epsilon _{R}$ which depends
on the root system and scales linearly with the inner product.
\begin{defn}\label{defn: e}
Let $n_{i}$ be the $i$th coefficient of the largest root
\[
\tilde{\alpha }=\sum _{i=1}^{n}n_{i}\alpha _{i}.
\]
We define
\[
\epsilon _{R}=\frac{1}{2}\LangRang{\tilde{\alpha },\tilde{\alpha }
} +\frac{1}{2}\sum _{i=1}^{n}n_{i}^{2}\LangRang{\alpha _{i},\alpha
_{i}}.
\]
\end{defn}
Note that in the case where $R$ is as in section~\ref{sec: intro},
namely of ADE type and the roots have norm square 2, then $\epsilon
_{R}=1+\sum _{i=1}^{n}n_{i}^{2}$ and we have that
\[
\epsilon _{R}=|G|
\]
where $G$ is the corresponding finite subgroup of $SU(2)$. This is a
consequence of the McKay correspondence, part of which implies that
$1, n_{1},\dotsc ,n_{n}$ are the dimensions of the irreducible
representations of $G$ (see
\cite[page~411]{Gonzalez-Sprinberg-Verdier}).

\subsection{The algebra $QH_{R}$}
Let 
\[
H_{R} = \znums \oplus \znums \alpha _{1}\oplus \dotsb \oplus \znums \alpha _{n}
\]
be the affine root lattice and let $QH_{R}$ be the free module over
$\znums [t][[q_{1},\dotsc ,q_{n}]]$ generated by $1,\alpha _{1},\dotsc
,\alpha _{n}$,
\[
QH_{R} = H_{R}\otimes \znums [t][[q_{1},\dotsc ,q_{n}]].
\]
We extend the pairing $\LangRang{\,,\,}$ to a $\qnums
[t,t^{-1}][[q_{1},\dotsc ,a_{n}]]$ valued pairing on $QH_{R}$ by
making 1 orthogonal to $\alpha _{i}$ and setting
\[
\LangRang{1,1} = \frac{-1}{t^{2}\epsilon _{R}}.
\]
For $\beta =\sum _{i=1}^{n}b_{i}\alpha _{i}$, we use the notation
\[
q^{\beta } = \prod _{i=1}^{n}q_{i}^{b_{i}}.
\]

\begin{thm}\label{thm: QH for any R}
Define a product operation $\star $ on $QH_{R}$ by letting 1 be the
identity and defining
\[
\alpha _{i}\star \alpha _{j} = -t^{2}\epsilon _{R} \LangRang{\alpha
_{i},\alpha _{j}} + \sum _{\beta \in R^{+}}\LangRang{\alpha _{i},\beta
}\LangRang{\alpha _{j},\beta ^{\vee }} t\frac{1+q^{\beta }}{1-q^{\beta
}}\beta .
\]
Then the product is associative, and moreover, it satisfies the
Frobenius condition
\[
\LangRang{x\star y,z} = \LangRang{x,y\star z}
\]
making $QH_{R}$ into a Frobenius algebra over the ring $\qnums
[t,t^{-1}][[q_{1},\dotsc ,q_{n}]]$.
\end{thm}

\smallskip

\begin{cor}\label{cor: W acts on QH by automorphism}
The Weyl group acts on $QH_{R}$ (and thus on $QH_{\cnums ^{*}}^{*}
(Y)$) by automorphisms. Namely, if we define
\[
g (q^{\beta })=q^{g\beta }
\]
for $g\in W$,
then for $v,w\in QH_{R}$ we have
\[
g (v\star w) = (gv)\star (gw).
\]
\end{cor}
\textsc{Proof:} Let $s_{k}$ be the reflection about the hyperplane
orthogonal to $\alpha _{k}$. By \cite[VI.1.6 Corollary~1]{Bourbaki},
$s_{k}$ permutes the positive roots other than $\alpha _{k}$. And
since the terms
\[
\frac{1+q^{\beta }}{1-q^{\beta }}\beta \quad \text{and}\quad
\LangRang{\alpha _{i},\beta }\LangRang{\alpha _{j},\beta^{\vee } }
\]
remain unchanged under $\beta \mapsto -\beta $, the effect of applying
$s_{k}$ to the formula for $\alpha _{i}\star \alpha _{j}$ is to
permute the order of the sum:
\begin{align*}
s_{k} (\alpha _{i}\star \alpha _{j})&= -t^{2} \epsilon _{R}
 \LangRang{\alpha _{i},\alpha _{j}} +\sum _{\beta \in R^{+}}
\LangRang{\alpha _{i},\beta } \LangRang{\alpha _{j},\beta ^{\vee
}}t\frac{1+q^{s_{k}\beta}}{1-q^{s_{k}\beta }}s_{k}\beta \\
&=-t^{2}\epsilon _{R}\LangRang{s_{k}\alpha _{i},s_{k}\alpha _{j}} + \sum
_{\beta \in R^{+}} \LangRang{\alpha _{i},s_{k}\beta } \LangRang{\alpha
_{j},s_{k}\beta ^{\vee }}t\frac{1+q^{\beta }}{1-q^{\beta }}\beta \\
&=s_{k} (\alpha _{i})\star s_{k} (\alpha _{j})
\end{align*}
and the Corollary follows.\qed 

\subsection{The proof of Theorem~\ref{thm: QH for any R}} When $R$ is
of ADE type and the pairing is normalized so that the roots have a
norm square of 2, then $QH_{R}$ coincides with $QH^{*}_{\cnums ^{*}}
(Y)$ and so Theorem~\ref{thm: QH for any R} for this case then follows
from Theorem~\ref{thm: product in QH}.

For any $R$, the Frobenius condition follows immediately from the
formulas for $\star $ and $\LangRang{\,,\,}$.

So what needs to be established in general is the associativity of the
$\star $ product. This is equivalent to the expression
\[
Ass^{R}_{xyuv} = \frac{1}{t^{2}}\LangRang{(x\star y)\star u,v}
\]
being fully symmetric in $\{x,y,u,v \}$. Written out, we have
\begin{multline*}
Ass^{R}_{xyuv} = -\epsilon _{R}\LangRang{x,y}\LangRang{u,v} \\
 + \sum _{\beta, \gamma \in R^{+} } \LangRang{x,\beta }
\LangRang{y,\beta } \LangRang{u,\gamma } \LangRang{v,\gamma }
\left(\frac{1+q^{\beta }}{1-q^{\beta }} \right)
\left(\frac{1+q^{\gamma }}{1-q^{\gamma }} \right)\LangRang{\beta
^{\vee },\gamma ^{\vee }}.
\end{multline*}
Recalling that $\epsilon _{R}$ scales linearly with the pairing, we
see that if $Ass^{R}_{xyuv}$ is fully symmetric in $\{x,y,u,v \}$,
then it remains so for any rescaling of the pairing.

To prove the associativity of $QH_{R}$ for root systems not of ADE
type, we reduce the non-simply laced case to the simply laced case.

Let $\{R, V,\LangRang{\,,\,} \}$ be an ADE root system and let $\Phi $
be a group of automorphisms of the Dynkin diagram. We construct a new
root system $\{R_{\Phi }, V^{\Phi },\LangRang{\,,\,}_{\Phi } \}$ as
follows. A somewhat similar construction can be found in
\cite[Section~10.3.1]{Springer}. Let
\[
V^{\Phi }\subset V
\]
be the $\Phi $ invariant subspace equipped with $\LangRang{\,,\,}_{\Phi
}$, the restriction of $\LangRang{\,,\,}$ to $V^{\Phi }$, and let the
roots of $R_{\Phi }$ be the $\Phi $ averages of the roots of $R$:
\[
R_{\Phi } = \left\{\overline{\alpha } = \frac{1}{|\Phi |}\sum _{g\in
\Phi } g\alpha , \alpha \in R \right\}.
\]
Then it is easily checked that $\{R_{\Phi },V^{\Phi
},\LangRang{\,,\,}_{\Phi } \}$ is an irreducible root system,
specifically of type given in the table:
\begin{center}
\begin{tabular}{|c|c|c|}
\hline $R$ &$\Phi $ &	$R_{\Phi }$\\
\hline $A_{2n-1}$&	$\znums _{2}$&	 $C_{n}$\\
$D_{n+1}$&	$\znums _{2}$&	 $B_{n}$\\
$E_{6}$&	$\znums _{2}$&	 $F_{4}$\\
$D_{4}$&	$\znums _{3}$&	 $G_{2}$\\
\hline 
\end{tabular}
\end{center}
Thus all the irreducible, reduced root systems arise in this way.

We will frequently use the fact that if $y\in V^{\Phi }$, then 
\begin{equation}\label{eqn: <x,y>=<xbar,y>}
\LangRang{x,y}=\LangRang{\overline{x},y}
\end{equation}
which easily follows from
$\LangRang{x,y}=\LangRang{gx,gy}=\LangRang{gx,y}$ for $g\in \Phi $.

We will also have need of the following two lemmas which we will prove
at the end of the section.
\begin{lemma}\label{lem: e for R = e for Rphi}
The constants defined in Definition~\ref{defn: e} coincide for the
root systems $R$ and $R_{\Phi }$:
\[
\epsilon _{R_{\Phi }} = \epsilon _{R}.
\]
\end{lemma}

\begin{lemma}\label{lem: sum over phi orbit of b is b dual}
Let $\beta \in R^{+}$ and let $\Phi \beta $ be the $\Phi $ orbit of
$\beta $. Then
\[
\sum _{\beta '\in \Phi \beta }\beta ' = \overline{\beta }^{\vee }.
\]
\end{lemma}

The simple roots of $R_{\Phi } $ are given by $\overline{\alpha
}_{i}$, the averages of the simple roots of $R$. Thus if 
\[
I=\{1,\dotsc ,n \}
\]
is the index set for the simple roots of $R$, then $\Phi $ acts on $I$ and 
\[
J=I/\Phi 
\]
is the natural index set for the simple roots of $R_{\Phi }$. For
$[i]\in J$, we let $\abar _{[i]}\in R_{\Phi }$ denote the simple root
given by $\abar _{i}$.

We specialize the variables $\{q_{i} \}_{i\in I}$ to variables
$\{\qbar_{[i]} \}_{[i]\in J}$ by setting
\begin{equation}\label{eqn: q specialization}
q_{i} = \qbar _{[i]}
\end{equation}
and it is straightforward to see that under the above specialization, 
\[
q^{\beta } = \qbar ^{\bbar }.
\]

Now let $R$ be an ADE root system whose roots have norm square 2. Then
$Ass^{R}_{xyuv}$ is fully symmetric in $\{x,y,u,v \}$. We specialize
the $q$ variables to the $\qbar $ variables as in equation~(\ref{eqn:
q specialization}) and we assume that $x,y,v,u\in V^{\Phi }$. Then
\begin{align*}
&Ass^{R}_{xyuv}+\epsilon _{R}\LangRang{x,y}\LangRang{u,v} =\\
&\quad \sum _{\beta ,\gamma
\in R^{+}} \LangRang{x,\beta }\LangRang{y,\beta }\LangRang{u,\gamma
}\LangRang{v,\gamma } \left(\frac{1+q^{\beta }}{1-q^{\beta }} \right)\left(\frac{1+q^{\gamma  }}{1-q^{\gamma  }} \right)\LangRang{\beta ^{\vee },\gamma ^{\vee }}\\
&=\sum _{\beta, \gamma \in
R^{+} } \LangRang{x,\bbar }\LangRang{y,\bbar }\LangRang{u,\gbar
}\LangRang{v,\gbar } \left(\frac{1+\qbar ^{\bbar }}{1-\qbar ^{\bbar }}
\right) \left(\frac{1+\qbar ^{\gbar }}{1-\qbar ^{\gbar }} \right)
\LangRang{\beta ,\gamma }\\
&=\sum _{\bbar ,\gbar  \in
R_{\Phi }^{+} } \LangRang{x,\bbar }\LangRang{y,\bbar }\LangRang{u,\gbar
}\LangRang{v,\gbar } \left(\frac{1+\qbar ^{\bbar }}{1-\qbar ^{\bbar }}
\right) \left(\frac{1+\qbar ^{\gbar }}{1-\qbar ^{\gbar }} \right)
\LangRang{\sum _{\beta '\in \Phi \beta }\beta' ,\sum _{\gamma '\in \Phi \gamma }\gamma' }\\
&=\sum _{\bbar, \gbar  \in
R_{\Phi }^{+} } \LangRang{x,\bbar }_{\Phi }\LangRang{y,\bbar }_{\Phi }\LangRang{u,\gbar
}_{\Phi }\LangRang{v,\gbar }_{\Phi } \left(\frac{1+\qbar ^{\bbar }}{1-\qbar ^{\bbar }}
\right) \left(\frac{1+\qbar ^{\gbar }}{1-\qbar ^{\gbar }} \right)
\LangRang{\bbar ^{\vee },\gbar ^{\vee }}_{\Phi }\\
&=Ass^{R_{\Phi }}_{xyuv} +\epsilon _{R_{\Phi }}\LangRang{x,y}_{\Phi }\LangRang{u,v}_{\Phi }
\end{align*}
and thus 
\[
Ass^{R_{\Phi }}_{xyuv} = Ass^{R}_{xyuv}
\]
is fully symmetric
in $\{x,y,u,v \}$ and the theorem is proved once we establish
Lemma~\ref{lem: e for R = e for Rphi} and Lemma~\ref{lem: sum over phi
orbit of b is b dual}.

\subsection{Proofs of Lemma~\ref{lem: e for R = e for Rphi} and Lemma~\ref{lem: sum over phi orbit of b is b dual}.}

We prove Lemma~\ref{lem: sum over phi orbit of b is b dual} first. If
$\beta $ is fixed by $\Phi $, the lemma is immediate. We claim that if
$\beta $ is not fixed then $\LangRang{\beta ,g\beta }=0$ for
non-trivial $g\in \Phi $. For simple roots, this follows from
inspection of the Dynkin diagrams and automorphisms which occur in the
table: a node is never adjacent to a node in its orbit. For other roots
this can also be seen from a direct inspection of the positive roots
(listed, for example, in \cite[Plates~I,IV--VII]{Bourbaki}).  For
$\beta $ not fixed by $\Phi $ we then have:
\begin{align*}
\LangRang{\bbar ,\bbar } &=\frac{1}{|\Phi |^{2}} 
\LangRang{\sum _{g}g\beta ,\sum _{h}h\beta }\\
&=\frac{1}{|\Phi |^{2}}\sum _{g}\LangRang{g\beta ,g\beta }\\
&=\frac{2}{|\Phi |}
\end{align*}
and Lemma~\ref{lem: sum over phi orbit of b is b dual} follows.

Note that the above formula generalizes to all roots $\bbar $ by
\[
\LangRang{\bbar ,\bbar } =2 \frac{\stab (\beta )}{|\Phi |}
\]
where $\stab (\beta )$ is the order of the stabilizer of the action of
$\Phi $ on $\beta $.

To prove Lemma~\ref{lem: e for R = e for Rphi} we must find the
coefficients of the longest root of $R_{\Phi }$. Since the longest
root of $R$ is unique, it is fixed by $\Phi $ and so it coincides with
the longest root of $R_{\Phi }$:
\begin{align*}
\overline{\tilde{\alpha }} =\tilde{\alpha } &= \sum _{i\in I}n_{i}\alpha _{i}\\
&=\sum _{[i]\in J} n_{[i]} \sum _{i'\in \Phi i} \alpha _{i}\\
&=\sum _{[i]\in J} n_{[i]} \abar _{[i]}^{\vee }\\
&=\sum _{[i]\in J} \frac{2n_{[i]}}{\LangRang{\abar _{[i]},\abar _{[i]}}} \abar _{[i]}.
\end{align*}
Thus we have 
\begin{align*}
2\epsilon _{R_{\Phi }} &=\LangRang{\overline{\tilde{\alpha }},\overline{\tilde{\alpha }}} +\sum _{[i]\in J} \left(\frac{2 n_{[i]}}{\LangRang{\abar _{[i]},\abar _{[i]}}} \right)^{2} \LangRang{\abar _{[i]},\abar _{[i]}}\\
&=\LangRang{\tilde{\alpha },\tilde{\alpha }} +\sum _{i\in I} \frac{\stab (\alpha _{i})}{|\Phi | } \frac{4n_{i}^{2}}{\LangRang{\abar _{i},\abar _{i}}}\\
&= \LangRang{\tilde{\alpha },\tilde{\alpha }} + \sum _{i\in I} 2 n_{i}^{2}\\
&=2\epsilon _{R}
\end{align*}
and Lemma~\ref{lem: e for R = e for Rphi} is proved.\qed 

\subsection{The root theoretic formula for triple
intersections. }\label{subsec: proof of prop}

Here we prove the root theoretic result required to finish the proof
of equation~\eqref{eqn: <aiajak>}.  Recall that $R$ is a root system
of ADE type normalized so that the roots have norm square 2.  We write
\[
g_{ij}=\LangRang{\alpha _{i},\alpha _{j}}.
\]

\begin{proposition}\label{prop: properties of Gijk} Let
\[
G_{ijk} = -\sum _{\beta \in R^{+}}\LangRang{\alpha _{i},\beta }
\LangRang{\alpha _{j},\beta } \LangRang{\alpha _{k},\beta }
\]
then $G_{ijk}$ satisfies the following properties.
\begin{enumerate}
\item [(i)] $G_{ijk} $ is symmetric in $\{i,j,k \}$,
\item [(ii)]$G_{ijk}$ is invariant under any permutation of indices
induced by a Dynkin diagram automorphism,
\item [(iii)]$G_{ijk}=0$ if $g_{jk}=0$,
\item [(iv)]$G_{ijk}=-8$ if $i=j=k$,
\item [(v)]$G_{iij}+G_{jji}=2$ if $g_{ij}=-1$,
\item [(vi)]$G_{ikk}+G_{jkk}=4$ if $i\neq j$ and $g_{ik}=g_{jk}=-1$.
\end{enumerate}
\end{proposition}

\textsc{Proof:} From the definition of $G_{ijk}$, properties (i) and
(ii) are clearly satisfied.  

Let $s_{k}$ be reflection about the hyperplane perpendicular to
$\alpha _{k}$ so that
\[
s_{k}\alpha _{i} = \alpha _{i} -g_{ik}\alpha _{k}.
\]
Since $s_{k}$ permutes the positive roots other than $\alpha _{k}$
\cite[VI.1.6 Corollary~1]{Bourbaki}, we get the following expression
for $G_{ijk}$:
\begin{align*}
G_{ijk} &= -2g_{ik}g_{jk}g_{kk}-\sum _{\beta \in R^{+}} \LangRang{\alpha _{i},s_{k}\beta } \LangRang{\alpha _{j},s_{k}\beta } \LangRang{\alpha _{k},s_{k}\beta }\\
&=-4g_{ik}g_{jk} - \sum _{\beta \in R^{+}} \LangRang{\alpha
_{i}-g_{ik}\alpha _{k},\beta } \LangRang{\alpha
_{j}-g_{jk}\alpha _{k},\beta } \LangRang{-\alpha _{k},\beta }\\
&=-4g_{ik}g_{jk} -G_{ijk} +g_{ik}G_{jkk} +g_{jk}G_{ikk} - g_{ik}g_{jk}G_{kkk}
\end{align*}
and so 
\begin{equation}\label{eqn: intermediate eqn for Gijk}
G_{ijk} = -2g_{ik}g_{jk}
+\frac{1}{2}\left(g_{ik}G_{jkk}+g_{jk}G_{ikk}-g_{ik}g_{jk}G_{kkk}
\right).
\end{equation}
Setting $i=j=k=n$ we obtain property (iv):
\[
G_{nnn}=-8
\]
which we can substitute back into equation~(\ref{eqn: intermediate eqn
for Gijk}) and then specialize $i=j=a$ to get
\begin{equation}\label{eqn: Gaak}
G_{aak} = 2g_{ak}^{2} + g_{ak}G_{akk}.
\end{equation}
Property (iii) then follows from equation~\eqref{eqn: intermediate eqn
for Gijk} and equation~\eqref{eqn: Gaak} and property (v) follows from
equation~\eqref{eqn: Gaak}.

For property (vi), observe that if $g_{ik}=g_{jk}=-1$ then $g_{ij}=0$
and so $G_{ijk}=0$ and equation~\eqref{eqn: intermediate eqn for Gijk}
then simplifies to prove property (vi).\qed

\section{Predictions for the orbifold invariants via the Crepant
Resolution Conjecture}\label{sec: CRC}

Let $G\subset SU (2)$ be a finite subgroup and let 
\[
\X =[\cnums ^{2}/G]
\]
be the orbifold quotient of $\cnums ^{2}$ by $G$. Recall that
\[
\pi :Y\to X
\]
is the minimal resolution of $X$, the singular variety underlying the
orbifold $\X $.

The Crepant Resolution Conjecture \cite{Bryan-Graber} asserts that
$F_{Y}$, the genus zero Gromov-Witten potential of $Y$, coincides with
$F_{\X }$, the genus zero orbifold Gromov-Witten potential of $\X $
after specializing the quantum parameters of $Y$ to certain roots of
unity and making a linear change of variables in the cohomological
parameters.

Using the Gromov-Witten computations of section~\ref{sec: GW theory of
Y}, we obtain a formula for $F_{Y}$. By making an educated guess for
the change of variables and roots of unity, and then applying the
conjecture, we obtain a prediction for the orbifold Gromov-Witten
potential of $\X $ (Conjecture~\ref{conj: prediction for F_X}). This
prediction has been verified in the cases where $G$ is $\znums _{2}$,
$\znums _{3}$, $\znums _{4}$ in
\cite{Bryan-Graber,Bryan-Graber-Pandharipande,Bryan-Jiang}
respectively, and recently it has been verified for all $\znums _{n}$
by Coates, Corti, Iritani, and Tseng \cite{CCIT-CRC}.

\subsection{The statement of the conjecture}
The variables of the potential function $F_{Y}$ are the quantum
parameters $\{q_{1},\dotsc ,q_{n} \}$ and cohomological parameters
$\{y_{0},\dotsc ,y_{n} \}$ corresponding the the generators
$\{1,\alpha _{1},\dotsc ,\alpha _{n} \}$ for $H_{\cnums ^{*}}^{*} (Y)$.

The potential function is the natural generating function for the genus
0 Gromov-Witten invariants of $Y$. It is defined by

\[
F_{Y} (q_{1},\dotsc ,q_{n},y_{0},\dotsc ,y_{n}) = \sum _{k_{0},\dotsc
,k_{n}} \sum _{A\in H_{2} (Y)} \LangRang{1^{k_{0}}\alpha
_{1}^{k_{1}}\dotsb \alpha _{n}^{k_{n}}}_{A}^{Y}
\,\,\frac{y_{0}^{k_{0}}}{k_{0}!}\dotsb \frac{y_{n}^{k_{n}}}{k_{n}!}
\,\, q^{A}.
\]

The potential function for the orbifold $\X =[\cnums ^{2}/G]$ depends
on variables $\{x_{0},\dotsc ,x_{n} \}$ which correspond to a basis
$\{1,\gamma _{1},\dotsc ,\gamma _{n} \}$ of $H_{\orb}^{*} (\X )$, the
orbifold cohomology of $\X $. The orbifold cohomology of $[\cnums
^{2}/G]$ has a natural basis which is indexed by conjugacy classes of
$G$. If $g\in G$ is an element of the group, we will write $x_{[g]}$
for the variable corresponding to the conjugacy class of $g$. There
are no curve classes in $\X $ and hence no quantum parameters so the
potential function is given by
\[
F_{\X } (x_{0},\dotsc ,x_{n}) = \sum _{k_{0},\dotsc
,k_{n}} \LangRang{1^{k_{0}}\gamma 
_{1}^{k_{1}}\dotsb \gamma _{n}^{k_{n}}}^{\X }
\,\,\frac{x_{0}^{k_{0}}}{k_{0}!}\dotsb \frac{x_{n}^{k_{n}}}{k_{n}!}.
\]
The conjecture states that there exists roots of unity $\omega
_{1},\dotsc \omega _{n}$ and an analytic continuation of $F_{Y}$ to
the points
\[
q_{i}=\omega _{i}
\]
such that the equality
\[
F_{Y} (\omega _{1},\dotsc ,\omega _{n},y_{0},\dotsc ,y_{n}) = F_{\X }
(x_{0},\dotsc ,x_{n})
\]
holds after making a (grading preserving) linear change of variables
\[
x_{i} = \sum _{j=0}^{n}L_{i}^{j}y_{j}.
\]

Thus to obtain a prediction for the potential $F_{\X }$, we must
determine the roots of unity $\omega _{i}$ and the change of variables
matrix\footnote{Note that our matrix $L$ here is the inverse of the
matrix called $L$ in \cite{Bryan-Graber}} $L$.

\subsection{The prediction.}
The only non-trivial invariants involving 1 are degree zero three
point invariants. We split up the potentials $F_{\X }$ and $F_{Y}$
into terms involving $x_{0}$ and $y_{0}$ respectively and terms
without $x_{0}$ and $y_{0}$ respectively.

Let $F^{0}_{Y}$ be the part of $F_{Y}$ with non-zero $y_{0}$ terms.
It follows from Lemma~\ref{lem: classical invs} that $F^{0}_{Y}$ is
given by
\[
F_{Y}^{0} = \frac{1}{t^{2}|G|}\frac{y_{0}^{3}}{3!} -
\frac{y_{0}}{2}\sum_{i,j=1}^{n}\LangRang{\alpha _{i},\alpha
_{j}}y_{i}y_{j}.
\]

Let $F^{0}_{\X }$ be the part of $F_{\X }$ with non-zero $x_{0}$ terms.
An easy localization computation shows that $F^{0}_{\X }$ is
given by
\[
F_{\X }^{0} = \frac{1}{t^{2}|G|}\frac{x_{0}^{3}}{3!} +
\frac{x_{0}}{2} \frac{1}{|G|}\sum _{g\in G,g\neq Id} x_{[g]}x_{[g^{-1}]}.
\]

Since the change of variables respects the grading, the terms in
$F_{Y}$ which are linear and cubic in $y_{0}$ must match up with the
terms in $F_{\X }$ which are linear and cubic in $x_{0}$. Consequently
we must have
\[
x_{0} = y_{0}
\]
and moreover, the change of variables must take the quadratic form
\begin{equation}\label{eqn: quadratic terms of F_X}
\frac{1}{|G|}\sum _{g\in G,g\neq Id} x_{[g]}x_{[g^{-1}]}
\end{equation}
to the quadratic form
\begin{equation}\label{eqn: quadratic terms of F_Y}
\sum _{i,j=1}^{n}-\LangRang{\alpha _{i},\alpha _{j}}y_{i}y_{j}.
\end{equation}

We can rewrite the above quadratic form in terms of the representation
theory of $G$ using the classical McKay correspondence
\cite{McKay} as follows. The simple roots $\alpha _{1},\dotsc ,\alpha _{n}$,
which correspond to nodes of the Dynkin diagram, also correspond to
non-trivial irreducible representations of $G$, and hence to their
characters $\chi _{1},\dotsc ,\chi _{n}$. Under this correspondence,
the Cartan paring can be expressed in terms of $\LangRang{\cdot |\cdot}$, the
natural pairing on the characters of $G$:
\begin{align*}
-\LangRang{\alpha _{i},\alpha _{j}} &= \LangRang{(\chi _{V}-2)\chi _{i}|\chi _{j}}\\
&=\frac{1}{|G|}\sum _{g\in G} (\chi _{V} (g)-2)\chi _{i}
(g)\overline{\chi }_{j} (g)
\end{align*}
where $V$ is the two dimensional representation induced by the
embedding $G\subset SU (2)$.

This discussion leads to an obvious candidate for the change of
variables. Namely, if we substitute
\begin{equation}\label{eqn: change of vars: x=Ly}
x_{[g]} = \sqrt{\chi _{V}-2} \,\,\sum _{i=1}^{n} \chi _{i} (g) y_{i}
\end{equation}
into equation~\eqref{eqn: quadratic terms of F_X} we obtain
equation~\eqref{eqn: quadratic terms of F_Y}. Since $\chi _{V} (g)$ is
always real and less than or equal to 2, we can fix the sign of the
square root by making it a positive multiple of $i$.

Thus we've seen that
\[
F_{Y}^{0} = F_{\X }^{0}
\]
under the change of variables given by equation~\eqref{eqn: change of
vars: x=Ly} and $x_{0}=y_{0}$. So from here on out, we set 
\[
x_{0}=y_{0}=0
\]
and deal with just the part of the potentials $F_{\X }$ and $F_{Y}$
not involving $x_{0}$ and $y_{0}$.

We apply the divisor axiom and the computations of section~\ref{sec:
GW theory of Y}:

\begin{align*}
F_{Y} &(y_{1},\dotsc ,y_{n},q_{1},\dotsc ,q_{n}) \\
&=\frac{1}{6}\sum
_{i,j,k=1}^{n} \LangRang{\alpha _{i}\alpha _{j}\alpha _{k}}_{0} y_{i}
y_{j}y_{k}+\sum _{A\in H_{2} (Y), A\neq 0} \LangRang{\,\,}_{A}q^{A}
e^{\sum _{i=1}^{n}y_{i}\int _{A}\alpha _{i}}\\
&= \frac{-t}{6} \sum _{i,j,k=1}^{n}\sum _{\beta \in R^{+}}
\LangRang{\alpha _{i},\beta } \LangRang{\alpha _{j},\beta }
\LangRang{\alpha _{k},\beta }y_{i}y_{j}y_{k}+\sum _{d=1}^{\infty }\sum
_{\beta \in R^{+}} \frac{2t}{d^{3}}q^{d\beta }e^{\sum
_{i=1}^{n}-d\LangRang{\alpha _{i},\beta }y_{i}}.
\end{align*}
Taking triple derivatives we get
\begin{align*}
\frac{\partial ^{3}F_{Y}}{\partial y_{i}\partial
y_{j}\partial y_{k}} =&-t\sum _{\beta \in R^{+}}\LangRang{\alpha
_{i},\beta } \LangRang{\alpha _{j},\beta } \LangRang{\alpha _{k},\beta
} \left(1+\frac{2q^{\beta }e^{{\sum _{i}-\LangRang{\beta ,\alpha
_{i}}y_{i}}}}{1-q^{\beta }e^{\sum _{i}-\LangRang{\beta ,\alpha
_{i}}y_{i}}} \right)\\
=&-t\sum _{\beta \in R^{+}} \LangRang{\alpha _{i},\beta }
\LangRang{\alpha _{j},\beta } \LangRang{\alpha _{k},\beta }
\frac{1+q^{\beta }e^{\sum _{i}-\LangRang{\beta ,\alpha
_{i}}y_{i}}}{1-q^{\beta }e^{\sum _{i}-\LangRang{\beta ,\alpha _{i}}y_{i}}}.
\end{align*}     
We specialize the quantum parameters to roots of unity by
\begin{equation}\label{eqn: roots of unity}
q_{j} = \exp\left(\frac{2\pi i n_{j}}{|G|} \right)
\end{equation}
where $n_{j}$ is the $j$th coefficient of the largest root as in
Definition~\ref{defn: e}. Note that $n_{j}$ is also the dimension of
the corresponding representation.

After specializing the quantum parameters, the triple derivatives of
the potential $F_{Y}$ can be expressed in terms of the function
\[
H (u) = \frac{1}{2i}\left(\frac{1+e^{i (u-\pi )}}{1-e^{i (u-\pi )}} \right) =
\frac{1}{2}\tan \left(\frac{-u}{2} \right)
\]
as follows:
\[
\frac{\partial ^{3}F_{Y}}{\partial y_{i}\partial
y_{j}\partial y_{k}} = -2i t\sum _{\beta \in R^{+}} \LangRang{\alpha
_{i},\beta } \LangRang{\alpha _{j},\beta } \LangRang{\alpha _{k},\beta
} H (Q_{\beta })
\]
where for $\beta = \sum _{j=1}^{n}b_{j}\alpha _{j}$ we define
\[
Q_{\beta } = \pi +\sum _{j=1}^{n}\left(\frac{2\pi n_{j}b_{j}}{|G|} + i
\LangRang{\beta ,\alpha _{j}}y_{j} \right).
\]
It then follows that 
\[
F_{Y} (y_{1},\dotsc ,y_{n}) =2t\sum _{\beta \in R^{+}} h (Q_{\beta })
\]
where $h (u)$ is a series satisfying
\[
h''' (u) = \frac{1}{2}\tan \left(\frac{-u}{2} \right).
\]

We can now make the change of variables given by equation~\eqref{eqn:
change of vars: x=Ly}. 
\begin{align*}
\sum _{j=1}^{n}i\LangRang{\beta ,\alpha _{j}} y_{j} &= \sum
_{j,k=1}^{n} ib_{k} \LangRang{\alpha _{k},\alpha _{j}} y_{j}\\
&= \sum _{j,k=1}^{n}  \frac{-i b_{k}}{|G|}\sum _{g\in G} (\chi _{V}
(g)-2)\,\overline{\chi}_{k} (g)\chi _{j} (g)y_{j}\\
&=\sum _{k=1}^{n}\frac{b_{k}}{|G|}\sum _{g\in G}\sqrt{2-\chi _{V}
(g)}\,\,\, \overline{\chi }_{k} (g) x_{[g]}.
\end{align*}
Substituting this back into $Q_{\beta }$ we arrive at our conjectural
formula for $F_{\X }$.

\begin{conjecture}\label{conj: prediction for F_X}
Let $F_{\X } (x_{1},\dotsc ,x_{n})$ denote the $\cnums ^{*}$
equivariant genus zero orbifold Gromov-Witten potential of the
orbifold $\X =[\cnums ^{2}/G]$ where we have set the unit parameter
$x_{0}$ equal to zero. Let $R$ be the root system associated to $G$ as
in section~\ref{sec: intro}. Then
\[
F_{\X } (x_{1},\dotsc ,x_{n}) = 2t \sum _{\beta \in R^{+}} h (Q_{\beta })
\]
where $h (u)$ is a series with
\[
h''' (u) = \frac{1}{2}\tan \left(\frac{-u}{2} \right)
\]
and 
\[
Q_{\beta } = \pi +\sum _{k=1}^{n}\frac{b_{k}}{|G|}\left(2\pi n_{k} +
\sum _{g\in G}\sqrt{2-\chi _{V} (g)}\,\, \overline{\chi }_{k} (g)
x_{[g]} \right)
\]
where $b_{k}$ are the coefficients of $\beta \in R^{+}$, $n_{k}$ are
the coefficients of the largest root, and $V$ is the two dimensional
representation induced by the embedding $G\subset SU (2)$.
\end{conjecture}
Note that the index set $\{1,\dotsc ,n \}$ in the above formula
corresponds to
\begin{enumerate}
\item simple roots of $R$,
\item non-trivial irreducible representations of $G$, and
\item non-trivial conjugacy classes of $G$.
\end{enumerate}
The index of a conjugacy class containing a group element $g$ is
denoted by $[g]$. Finally note that the terms of degree less than
three are ill-defined for both the potential $F_{\X }$ and our
conjectural formula for it.

The above conjecture has been proved in the cases where $G$ is $\znums
_{2}$, $\znums _{3}$, $\znums _{4}$ in
\cite{Bryan-Graber,Bryan-Graber-Pandharipande,Bryan-Jiang}
respectively, and recently it has been verified for all $\znums _{n}$
by Coates, Corti, Iritani, and Tseng \cite{CCIT-CRC}.

We have also performed a number of checks of the conjecture for
non-Abelian $G$. Many of the orbifold invariants must vanish by monodromy
considerations, and our conjecture is consistent with this
vanishing. One can geometrically derive a relationship between
some of the orbifold invariants of $[\cnums ^{2}/G]$ and certain
combinations of the orbifold invariants of $[\cnums ^{2}/H]$ when $H$
is a normal subgroup of $G$. This leads to a simple relationship
between the corresponding potential functions which we have checked is
consistent with our conjecture.

\bibliography{mainbiblio} \bibliographystyle{plain}

\begin{thebibliography}{10}

\bibitem{Be-Fa}
K.~Behrend and B.~Fantechi.
\newblock The intrinsic normal cone.
\newblock {\em Invent. Math.}, 128(1):45--88, 1997.

\bibitem{Bertram}
Aaron Bertram.
\newblock Another way to enumerate rational curves with torus actions.
\newblock {\em Invent. Math.}, 142(3):487--512, 2000.

\bibitem{Bourbaki}
N.~Bourbaki.
\newblock {\em \'{E}l\'ements de math\'ematique. {F}asc. {XXXIV}. {G}roupes et
  alg\`ebres de {L}ie. {C}hapitre {IV}: {G}roupes de {C}oxeter et syst\`emes de
  {T}its. {C}hapitre {V}: {G}roupes engendr\'es par des r\'eflexions.
  {C}hapitre {VI}: syst\`emes de racines}.
\newblock Actualit\'es Scientifiques et Industrielles, No. 1337. Hermann,
  Paris, 1968.

\bibitem{Bryan-Gholampour3}
Jim Bryan and Amin Gholampour.
\newblock {The quantum McKay corresondence}.
\newblock In preparation.

\bibitem{Bryan-Graber}
Jim Bryan and Tom Graber.
\newblock The crepant resolution conjecture.
\newblock arXiv: math.AG/0610129.

\bibitem{Bryan-Graber-Pandharipande}
Jim Bryan, Tom Graber, and Rahul Pandharipande.
\newblock The orbifold quantum cohomology of {$\mathbf{C^2/Z_3}$ and Hurwitz
  Hodge integrals}.
\newblock arXiv:math.AG/0510335, to appear in Journal of Alg. Geom.

\bibitem{Bryan-Jiang}
Jim Bryan and Yunfeng Jiang.
\newblock {The Crepant Resolution Conjecture for the orbifold
  $\mathbf{C^2/Z_4}$}.
\newblock In preparation.

\bibitem{BKL}
Jim Bryan, Sheldon Katz, and Naichung~Conan Leung.
\newblock Multiple covers and the integrality conjecture for rational curves in
  {C}alabi-{Y}au threefolds.
\newblock {\em J. Algebraic Geom.}, 10(3):549--568, 2001.
\newblock Preprint version: math.AG/9911056.

\bibitem{CCIT-CRC}
Tom Coates, Alessio Corti, Hiroshi Iritani, and Hsian-Hua Tseng.
\newblock {The Crepant Resolution Conjecture for Type A Surface Singularities}.
\newblock arXiv:0704.2034v1 [math.AG].

\bibitem{Gonzalez-Sprinberg-Verdier}
G.~Gonzalez-Sprinberg and J.-L. Verdier.
\newblock Construction g\'eom\'etrique de la correspondance de {M}c{K}ay.
\newblock {\em Ann. Sci. \'Ecole Norm. Sup. (4)}, 16(3):409--449 (1984), 1983.

\bibitem{Ka-Mo}
Sheldon Katz and David~R. Morrison.
\newblock Gorenstein threefold singularities with small resolutions via
  invariant theory for {W}eyl groups.
\newblock {\em J. Algebraic Geom.}, 1(3):449--530, 1992.

\bibitem{Maulik-An}
Davesh Maulik.
\newblock {Gromov-Witten Theory of A-resolutions}.
\newblock In preparation.

\bibitem{McKay}
John McKay.
\newblock Graphs, singularities, and finite groups.
\newblock In {\em The Santa Cruz Conference on Finite Groups (Univ. California,
  Santa Cruz, Calif., 1979)}, volume~37 of {\em Proc. Sympos. Pure Math.},
  pages 183--186. Amer. Math. Soc., Providence, R.I., 1980.

\bibitem{Reid-asterisque}
Miles Reid.
\newblock La correspondance de {M}c{K}ay.
\newblock {\em Ast\'erisque}, (276):53--72, 2002.
\newblock S\'eminaire Bourbaki, Vol.\ 1999/2000.

\bibitem{Springer}
T.~A. Springer.
\newblock {\em Linear algebraic groups}, volume~9 of {\em Progress in
  Mathematics}.
\newblock Birkh\"auser Boston Inc., Boston, MA, second edition, 1998.

\end{thebibliography}

\end{document}